\documentclass[aps,amssymb,pre,groupedaddress,twocolumn]{revtex4-1}

\usepackage{graphicx}
\usepackage{dcolumn}
\usepackage{bm}
\usepackage{amsfonts}
\usepackage{amsmath, amssymb}
\usepackage{color}
\usepackage{multirow}
\usepackage{longtable}

\newcommand\qed{\phantom{\underline{y}}\hfill\hfill$\square$}

\newcommand\1{{\bf 1}}

\newcommand\E{{\mathbb E}}

\newcommand{\weg}[1]{}
\newcommand\ed{\stackrel{d}{=}}

\begin{document}

\title{Consistency of detrended fluctuation analysis}
\author{O. L{\o}vsletten}
\address{Department of Mathematics and Statistics, UiT Artic University of Norway, Norway.}

\begin{abstract}
The scaling function $F(s)$ in detrended fluctuation analysis (DFA) scales as $F(s)\sim s^{H}$ for stochastic processes with Hurst exponent $H$. This scaling law is proven for  stationary stochastic processes with $0<H<1$, and non-stationary stochastic processes with $1<H<2$.  
For $H<0.5$ it is observed that  the asymptotic (power-law) auto-correlation function (ACF) scales as $\sim s^{1/2}$.
It is also demonstrated that the fluctuation function in DFA is equal in expectation to: i) a weighted sum of the ACF,  ii) a weighted sum of the second order structure function. These results enable us to compute the exact finite-size bias for signals that are scaling, and to employ  DFA in a meaningful sense for signals that do not exhibit power-law statistics. The usefulness is illustrated by  examples where it is demonstrated that a previous suggested modified DFA will increase the bias for signals with Hurst exponents $H>1$. As a final application of the new developments,  a new estimator $\hat F(s)$  is proposed. This estimator  can handle missing data in regularly sampled time series without the need of interpolation schemes. Under mild regularity conditions, $\hat F(s)$ is equal in expectation to the fluctuation function $F(s)$ in the gap-free case.

\end{abstract}
\maketitle
\section{Introduction}
Detrended fluctuation analysis (DFA) was introduced in a study of long-range dependence in DNA sequences \cite{Peng:1994be}. It has later been applied in a wide range of scientific disciplines \citep{Kantelhardt:2001dr}. Some recent examples are found in scientific studies of climate \citep{Lovsletten:2016bja}, finance \citep{Lahmiri:2015ch} and medicine \citep{Chiang:2016el}. The most common usage of DFA is to estimate the Hurst exponent. The assumption is then that the second moment of the fluctuations of the input signal, after these have been averaged over a time scale $s$, is a power-law function of $s$. This property is called {\em scale invariance}, or just {\em scaling}. In the context of DFA, the scaling assumption implies that  the DFA fluctuation function $F(s)$ takes the form of  a power-law 
\begin{equation}
\E\, F^2(s)\sim s^{2H},
 \label{mr}
\end{equation}
where $\E $ denotes the expectation, i.e., the ensemble mean.

Important examples of stochastic processes $X(t)$ with scaling properties are self-similar and multifractal models, see e.g., \citep{Lovsletten:2012ga}. For this large class of models, the existing $q$-moments
satisfy $\E \left|X(t+t_0)-X(t_0)\right|^q\propto t^{\zeta(q)}$. In particular, if the variance is finite, the second moments are scaling and  the Hurst exponent $H$ is defined by the relation $\zeta(2)=2H-2$. The power-law of the DFA fluctuation function in this case ($1<H<2$) has been established empirically. A mathematical proof has not been published prior to this paper, except for random walks ($H=1.5$) \citep{Holl:2015it,Holl:2015gv}.

For stationary stochastic processes $X(t)$ with scaling  second-moments, the Hurst exponent is in the range $0<H<1$. For $H=0.5$, $X(t)$ is white noise, while $H\neq 1/2$ implies an auto-correlation function (ACF) $\rho(\tau)$ on the form
\citep{Beran:1994uu}

\begin{equation}
\rho(\tau)\sim H(2H-1)\tau^{2H-2}.
\end{equation}
For $H<1/2$ the ACF is negative for all time lags $\tau \neq 0$, while for $H>1/2$ the ACF is positive. Moreover, in the persistent case ($H>1/2$), the ACF decays so slowly that the series $\sum_{\tau=-\infty}^{\infty} {\rho(\tau)}$ diverges.


In the case of a stationary input signal $X(t)$, with Hurst exponent $0<H<1$, Eq.~\eqref{mr} has been partly proven in the past. \citet{TAQQU:1995dt} constructed a proof for DFA1. DFA$m$, or DFA of order $m$, means that a $m$'th order polynomial is applied in the DFA algorithm (Section \ref{sec3}A).
For Hurst exponents restricted to the range $0.5<H<1$, the proof has been extended to include higher-order polynomials $m\geq 1$ \cite{Holl:2015gv}. In this paper a new observation is made:  for $0<H<0.5$, in order for Eq.~\eqref{mr} to be satisfied, only the {\em exact} auto-covariance function (acvf)  gives the correct result. If  the asymptotic acvf is employed, then $\E\, F^2(s)\sim s$.

For stationary signals, \citet{Holl:2015gv} showed that the squared DFA fluctuation function is equal in expectation to a weighted sum of the acvf $\gamma(\cdot)$: 
\begin{equation}
\mathbb E F^2 (s)= \gamma(0) G(0,s)s^{-1}+2 s^{-1}\sum_{j=1}^{s-1} G(j,s)\gamma(j), 
\label{EACF}
\end{equation}
where the weight function $G(j,s)$ will be defined in Section~\ref{sec3}. 
In this paper a more general result is presented;
\begin{equation}
\mathbb E F^2(s)= -\frac{1}{s} \sum_{j=1}^{s-1} G(j,s)S(j),
\label{ES}
\end{equation}
where $S(t)=\E [X(t+t_0)-X(t_0)]^2$, which also holds for non-stationary stochastic processes with stationary increments. The quantity $S(t)$ is known as the variogram. We note that the relationship between DFA and the power spectral density was derived, partly based on numerical calculations, by
\citet{Heneghan:2000iz}.

Eqs.~\eqref{EACF}~and~\eqref{ES} have applications beyond proving Eq.~\eqref{mr}.  For instance, one can compute the exact finite-size bias for scaling  signals, and  make meaningful use of DFA for signals that are not scaling. 
In \citet{Kantelhardt:2001dr}, the bias of DFA for stochastic processes with Hurst exponents in the range $0.5\leq H<1$ was found by means of Monte Carlo  simulations, using long time series of synthetically generated fractional Gaussian noises.  An analytical study of the behaviour of DFA for auto-regressive processes of order one (AR(1)) was investigated in \citet{Holl:2015gv}. In the present paper the usage of   Eqs.~\eqref{EACF}~and~\eqref{ES} is demonstrated  by simple extensions of the aforementioned examples.

An important application of the theoretical developments  presented in this paper is the construction of estimators (modifications of the DFA fluctuation function) that can handle missing data in regularly sampled time series. One simple way of handling missing data is to apply linear interpolation, random re-sampling, or mean filling. 
However, this will typically destroy, or add, artificial correlations to the time series under study. The effect on DFA using these three gap-filling techniques was examined in \citet{Wilson:2003ij} for signals with Hurst exponents $0<H<1$. It was found that these interpolation schemes introduce significant deviation from the expected scaling. In contrast, the modified fluctuation functions  proposed here have the property of equality in expectation to the fluctuation function in the gap-free case.  For the wavelet variance, estimators that can handle missing data in a proper statistical way was presented by  \citet{Mondal:2008cg}. These wavelet variances are similar in construction to the DFA estimators presented here.


This paper is organised as follows. In Section~\ref{sec2} the definition of Hurst exponent adopted in this paper is reviewed. Examples of stochastic processes with well-defined Hurst exponents are given. 
Section~\ref{sec3} presents the relationship between the DFA fluctuation function and the acvf and variogram, and the proof of Eq.~\eqref{mr}.
Examples of applications are given in Section~\ref{sec3}: Bias for scaling signals, DFA of Ornstein-Uhlenbeck processes, and modification of the DFA fluctuation function to handle missing data.



\section{Hurst exponent \label{sec2}}
\subsection{Definition and properties}
Let $X(t)$ be a stochastic process with mean $\mathbb E X(t)=0$. If
\\ i) $X(t)$ is non-stationary with stationary increments and
\begin{equation*}
\mathbb E \left[X(t+t_0)-X(t_0)\right]^2\propto t^{2 H-2},
\label{DefH1}
\end{equation*}
or \\ii)
 $X(t)$ is stationary and
\begin{equation*}
\mathbb E [Y(t+t_0)-Y(t_0)]^2\propto t^{2 H},\quad  Y(t)=\sum_{k=1}^t X(k) ,
\label{DefH2}
\end{equation*}
then we define $H$ to be the Hurst exponent of the process $X(t)$. The Hurst exponent determines the correlation at all time scales. Assume that $X(t)$ has Hurst exponent $1<H<2$, i.e., $X(t)$ is non-stationary. Then,
\begin{equation*}
2X(t)X(s)= X(t)^{2}+X(s)^{2}-\{X(t)-X(s)\}^{2}.
\end{equation*}
By stationary increments;
\begin{equation*}
\E\{X(t)-X(s)\}^{2} =\E X(|t-s|)^{2},
\end{equation*}
it follows that
\begin{equation}
\E X(t)X(s)=\frac{\sigma^2}{2}\left\{|s|^{2h}+|t|^{2h}-|t-s|^{2h}\right\},
\label{CovMotion}
\end{equation}
with $\E X(1)^2=\sigma^2$ and $h=H-1$. The increment process $\Delta X(t)=X(t)-X(t-1)$ has Hurst exponent $h$.
The acvf $\gamma(\tau)$ of the increments follows from \eqref{CovMotion}, and is given by
\begin{equation}
\gamma(\tau)= \frac{\sigma^2}{2}\left(|\tau + 1|^{2h} - 2  |\tau|^{2h} +|\tau - 1|^{2h}\right). 
\label{ACF}
\end{equation}
For $h=1/2$ the increment process is white noise, while for $h\neq 1/2$ the acvf is asymptotically a power-law;
\begin{equation*}
\gamma(\tau)
\sim   \frac{\sigma^2}{2}\frac{d^2}{d\tau^2} t^{2h}=\sigma^2 h(2h-1)\tau^{2h-2} ,
\end{equation*}
as  $\tau\rightarrow \infty$. Thus, $h \neq1/ 2$ implies dependent increments. Choosing $0<h<1/2$ results in negatively correlated increments, while for $h>1/2$ the increments are persistent. Moreover, in the persistent case, the acvf decays so slowly that the series $\sum_{\tau=-\infty}^{\infty} {\gamma(\tau)}$ diverges. 

\subsection{Examples}
Hurst exponent in the range $1<H<2$ and $X(t)$  Gaussian distributed  defines the class of fractional Brownian motions (fBm's). The corresponding increment process is known as a fractional Gaussian noise (fGn's) \citep{Mandelbrot:1968el}. An fBm is an example of a self-similar process.
By definition self-similar processes $X(t)$, with self-similar exponent $h$, satisfy the the scale-invariance
\begin{equation}
X(a t)\ed M(a) X(t),
\label{multifractal}
\end{equation}
where $M(a)=a^h$ \cite{embrechts2002selfsimilar}, and  ``$\ed$" denotes equality in finite-dimensional  distributions. The class of log-infinitely divisible multifractal processes \citep{Muzy:2002cq,Bacry:2003jr} also satisify Eq.~\eqref{multifractal}, but now $M(a)$ is random variable with an arbitrarily log-infinitely divisible distribution. The scaling law Eq.~\eqref{multifractal} implies $\E \left|X(t+t_0)-X(t_0)\right|^q\propto t^{\zeta(q)}$. Thus, if the second moments exist the Hurst exponent $H$ is given by the relation $\zeta(2)=2H-2$. These examples are summarized in Table~\ref{T1}. We emphasize that neither multifractality nor self-similarity is needed to have a process with well-defined Hurst exponent. An example is the class of smoothly truncated L\'evy flights (STLF's) \citep{Koponen:1995ev}. For STLF's all moments exist, and the property of stationary and independent increments implies a Hurst exponent $H=1.5$. The STLF behaves like a L\`evy flight on small time scales, while on long time scales, the statistics are close to Brownian motion \cite{Terdik:2006ch}. Thus, it is clearly neither self-similar nor multifractal, which was proven in \cite{Rypdal:2016bd}. 


 
\begin{table}
\begin{tabular}{| l| c | }
 \hline
 Stochastic process & Hurst exponent  \\
 \hline
  White noise &  $H=1/2$  \\
  Random walks & $H=3/2$  \\
 fractional Gaussian noise  & $0<H<1$\\
  fractional Brownian motion  & $1<H<2$\\
$h$-selfsimilar processes  & $H=h+1$\\ 
Scaling function $\zeta(q)$& $H=\zeta(2)/2+1$\\
 \hline
\end{tabular}
\caption{Examples of stochastic processes with well-defined Hurst-exponents $H$. Finite variance is assumed in all examples.  
}
\label{T1}
\end{table}

\section{Detrended fluctuation analysis \label{sec3}}
\subsection{DFA algorithm}
Let $X(1),X(2),\ldots,X(n)$ be the input to DFA. The first step in DFA is to construct the profile
\begin{equation*}
Y(t)=\sum_{k=1}^tX(k). 
\label{profile}
\end{equation*}
For a given scale $s$ one considers windows of length $s$. In each window a polynomial of degree $m$ is fitted to the profile. Subtracting the fitted polynomial from the profile gives a set of residuals. From these residuals the variance is computed. We denote by $F_t^2(s)$ the residual variance. The squared fluctuation function $F^2$ is the average of $F_t^2$. To express the residual variance mathematically, we introduce some notation. Define the vector $\mathbf{Y}(t)=[Y(t+1),Y(t+2),\ldots,Y(t+s)]^T$.
Let $B$ be the $(m+1)\times s$ design matrix in the ordinary least square (OLS) regression. 
That is, row $k$ of $B$ is the vector $(1^{k-1},2^{k-1},\ldots,s^{k-1})$. Define
\begin{equation}
Q=B^T\left(BB^T\right)^{-1} B, 
\label{hatmatrix}
\end{equation}
which is known as the hat matrix in statistics.
The residual variance is given by
\begin{equation}
F_t^2(s)=\frac{1}{s}\,\mathbf{Y}(t)^T(I-Q) \mathbf{Y}(t),
\label{DFAdef}
\end{equation}
where $I$ is the $(s\times s)$ identity matrix.
\subsection{How DFA relates to variogram and acvf}
\begin{figure}
\begin{center}
\includegraphics[width=\linewidth]{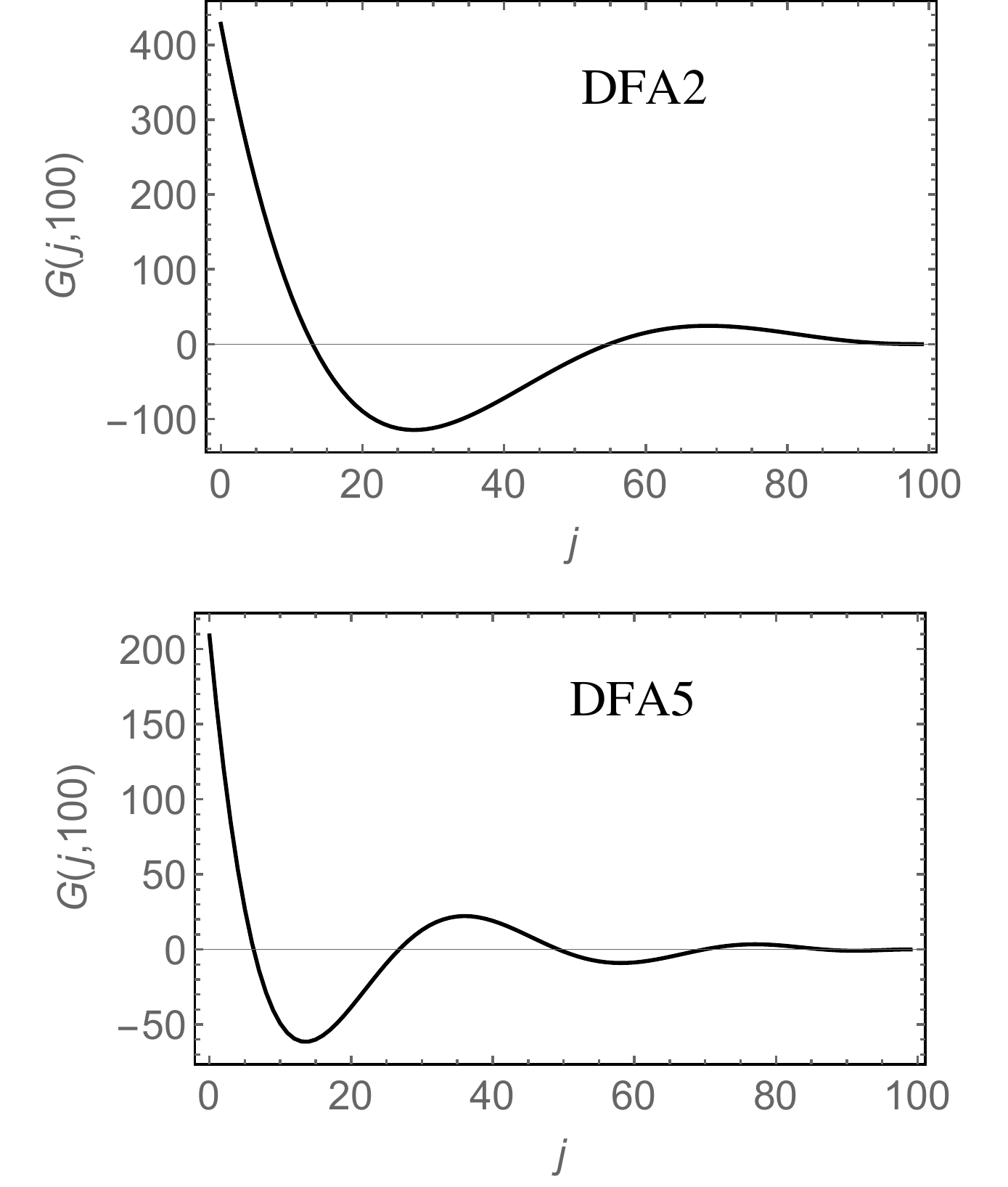}
\caption{Shows the map $j\mapsto G(j,100)$ for DFA2 (top figure) and DFA5 (bottom figure). The weight function $G(j,s)$ is defined in the text.   }
\label{fig1}
\end{center}
\end{figure} 

It is convenient to express the squared fluctuation function explicitly in terms of the input series.
 Let $\mathbf{X}(t)=[X(t+1),X(t+2),\ldots,X(t+s)]^T$. We define the $s\times s$ matrix $D$ by letting element $(i,j)$ of $D$ be  unity if $i\geq j$ and zero otherwise. Left-multiplying $D$ with $\mathbf{X}(t)$ gives the vector of cumulative sums $(X(t+1),X(t+1)+X(t+2),\ldots,\sum_{k=1}^sX(t+k))$. Define $$A=D^T(I-Q)D,$$
and let $a_{k,j}$ be element $(k,j)$ of the matrix $A$. 
The fluctuation function can be written
\begin{eqnarray}
F^2_t(s)&=&\frac{1}{s}\mathbf{X}(t)^TA\mathbf X (t)\nonumber\\
&=&\frac{1}{s}\sum_{k=1}^s\sum_{j=1}^s a_{k,j}X(t+k)X(t+j). 
 \label{F2}
\end{eqnarray}
In the definition of DFA the profile is constructed for the entire time series prior to windowing. Eq.~\eqref{F2} states that constructing the profile in each window gives the same residual variance (squared fluctuation function). 
Another form of the residual variance is
\begin{equation}
F^2_t(s)=-\frac{1}{2s}\sum_{k=1,j=1}^s a_{k,j} [X(t+k)-X(t+j)]^2.
\label{F3}
\end{equation}
The proofs of Eqs.~\eqref{F2} and \eqref{F3} can be found in Appendix~\ref{VA}. 

In the sequel we make the assumption that $X(t)=T(t)+Z(t)$, where $Z(t)$ is a stochastic process with mean zero and acvf $\gamma(t,s)$. The deterministic part $T(t)$ of the input signal is a trend modeled as a polynomial of order $q$ less than the order $m$ of DFA (see 
Appendix~\ref{VA} for precise form of the trends). For simplicity we consider in this section  the case $T(t)=0$, and postpone to Appendix~\ref{VA} to show that the results are valid also for trends with $q<m$.

By applying the expectation operator to Eq.~\eqref{F2}, it is seen that 
\begin{equation}
\mathbb E F_t^2(s)= \frac{1}{s}\sum_{k=1}^s\sum_{j=1}^s a_{k,j}\gamma(t+k,t+j).
 \label{acvf1}
\end{equation}

If we add the further restriction of stationarity of the process $X(t)$, Eq.~\eqref{acvf1} simplifies to
Eq.~\eqref{EACF}, with $\gamma(t)=\gamma(0,t)$ and 
\begin{equation}
G(j,s)=\sum_{k=1}^{s-j} a_{k,k+j}.
\label{G}
\end{equation}

While Eq.~\eqref{acvf1} appears to be  time-dependent when $X(t)$ is non-stationary with stationary increments, this is not the case. To establish that $\mathbb E F_t^2(s)$ does not depend on the window $t$, one can apply the expectation operator to Eq.~\eqref{F3}. The result is Eq.~\eqref{ES}, with $G(j,s)$ the weight function defined in  
Eq.~\eqref{G}. Since $\E F^2_t(s)$ does not depend on the window $t$, it follows that
 $\E F^2_t(s)= \E F^2(s)$.

The weight functions $G(j,s)$ can be computed exactly. In this work this has been done by means of Mathematica.
 The weight function for DFA1 and DFA2 are listed in Table~\ref{Tabell2}, while the map $j\mapsto G(j,100)$ for DFA2 and DFA5 are shown in Fig.~\ref{fig1}.

  \begin{longtable*}[!b]{c c}
\\
 
\hline
$m$ & $G(j,s)$ \\ 
\hline
\hline
1& $\frac{(j-s-1) (j-s) (j-s+1) \left(3 j^2+9 j s-2 s^2+8\right)}{30 s \left(s^2-1\right)}$ \\ 
 2& $-\frac{(j-s-1) (j-s) (j-s+1) \left(10 j^4+30 j^3 s+2 j^2 \left(9 s^2+19\right)+2 j s
   \left(67-13 s^2\right)+3 \left(s^4-13 s^2+36\right)\right)}{70 s \left(s^4-5
   s^2+4\right)}$ \\
   \hline \\
   
   \caption{The weight function $G(j,s)$ for DFA1 and DFA2. 
   }
   \label{Tabell2}
\end{longtable*}

  \begin{longtable*}{c cccccccccccccccccccccccccc}
\\
 \hline
 $q$& 0&1 & 2 & 3 & 4 & 5 & 6 & 7 & 8 & 9 & 10 & 11 & 12 & 13 & 14 & 15 \\
  \hline
  \hline
   \\[.01pt]
$m=1 $&  $\frac{1}{15}$ &$ -\frac{1}{2} $& $1$ &$ -\frac{2}{3}$ &$ 0$ & $\frac{1}{10}$ \\[2pt] 
$m=  2$ &$\frac{3}{70}$ &$ -\frac{1}{2}$ & $\frac{3}{2}$ &$ -\frac{3}{2}$ & $0$ &$ \frac{3}{5} $&$ 0 $&
   $-\frac{1}{7}$   \\[2pt] 
 $m= 3$ & $\frac{2}{63}$ & $-\frac{1}{2}$ &$ 2$ & $-\frac{8}{3} $&$ 0 $& $2$ & $0$ & $-\frac{8}{7}$ &$ 0$ &
  $ \frac{5}{18} $ \\[2pt] 
$m= 4$ & $\frac{5}{198} $&$ -\frac{1}{2} $&$ \frac{5}{2} $&$ -\frac{25}{6} $&$ 0$ & $5$ &$ 0$ & $-5$ &$ 0$ &
  $ \frac{25}{9} $& $0$ & $-\frac{7}{11}$  \\[2pt] 
$m=5$ &   $\frac{3}{143} $& $-\frac{1}{2}$ &$ 3$ &$ -6$ &$ 0$ & $\frac{21}{2}$ &$ 0$ &$ -16$ & $0$ & $15$ & $0$ &
  $ -\frac{84}{11}$ &$ 0$ &$ \frac{21}{13}$\\[2pt] 
$m=6$ &$\frac{7}{390}$ &$ -\frac{1}{2} $&$ \frac{7}{2}$ &$ -\frac{49}{6} $&$ 0$ &$ \frac{98}{5}$ & $0$ & $-42$
   & $0$ & $\frac{175}{3}$ & $0$ &$ -49 $&$ 0$ & $\frac{294}{13} $& $0$ &$ -\frac{22}{5}$ \\[1pt] \hline
\\
\caption{The coefficients $\{d_q\}$ for DFA of order $m= 1,2,\ldots,6$.}
\label{Tabell3}
  \end{longtable*}

  \subsection{Proof of DFA scaling}
We are now in a position to prove  \begin{equation}
\mathbb E F^2(s)\sim \lambda_{m,H} s^{2H},
\label{MR2}
\end{equation}
for input signals $X(t)$ with Hurst exponent $H\in\{(0,1)\cup (1,2) \}$. We assume $\E X(t)^2=1$. In the appendix we derive the asymptotic weight function $G_\text{asym}(j,s)\sim G(j,s)$, which takes the form
$$G_\text{asym}(j,s)
=
\begin{cases}
\sum_{q=0}^{2m+3} s^{2-q}j^qd_q& \text{if     } j>0, \\
d_0 s^2& \text{if     } j=0. \\
\end{cases}
$$
Expressions for the coefficients $\{d_q\}$ can be found in Appendix~\ref{VB}. 
The values of $\{d_q\}$ for orders $m\leq 6$ are listed in Table~\ref{Tabell3}. 

For $H>1$, using Eq.~\eqref{ES} and the asymptotic weight function, yields the asserted scaling:
\begin{eqnarray*}
\E F(s)^2&\sim& -s^{2H} \sum_{q=0}^{2m+3}s^{1-q} d_q \sum_{j=1}^{s-1}   j^{q+2H-2} 
\\&\sim& - s^{2H}\sum_{q=0}^{2m+3} \frac{d_q}{q+2H-1}. 
\end{eqnarray*}
In the stationary case $0<H<1$,  using Eq.~\eqref{EACF}, we find
\begin{eqnarray}
\E F^2(s) &\sim & s\, d_0 \gamma(0) +2 \sum_{q=0}^{2m+3}s^{1-q} d_q \sum_{j=1}^{s-1}   j^q \gamma(j)
\label{asymF}.
\end{eqnarray}
White noise ($H=1/2$) is trivial:
\begin{eqnarray*}
\E F(s)^2 &\sim & s\, d_0 .
\end{eqnarray*}
For $1/2<H<1$, in Eq.~\eqref{asymF} the second term dominates the first term. Coupled with the asymptotic form of the acvf we have 
\begin{eqnarray*}
\E F(s)^2&\sim&2 H(2H-1)\sum_{q=0}^{2m+3}s^{1-q} d_q \sum_{j=1}^{s-1}   j^{q+2H-2} \\
&\sim &s^{2H}2 H(2H-1)\sum_{q=0}^{2m+3} \frac{d_q}{q+2H-1}. \\
\end{eqnarray*}
For $0<H<1/2$ the leading terms in Eq.~\eqref{asymF} scale as $s$, 
but those leading terms cancel and we end up with the same formula as above. To see this, denote by $\rho(\tau)$ the auto-correlation function (ACF), i.e., $\rho(\tau)=\gamma(\tau)/\gamma(0)$. It is well-known that $\sum_{j=-\infty}^\infty \rho(\tau) =0$ (e.g., \citep{Beran:1994uu}). Since $\rho(0)=1$, and the ACF is a symmetric function, we have $-\gamma(0)/2=\sum_{j=1}^\infty \gamma(\tau)$. Thus 
\begin{eqnarray*}
\E F(s)^2&\sim& -2  d_0 s \sum_{j=s}^{\infty} \gamma(j) +2 \sum_{q=1}^{2m+3}s^{1-q} d_q \sum_{j=1}^{s-1}   j^q \gamma(j)\\
&\sim& -2  d_0 s H(2H-1) \sum_{j=s}^{\infty} j^{2H-2} \\
&&+2 H(2H-1)\sum_{q=1}^{2m+3}s^{1-q} d_q \sum_{j=1}^{s-1}   j^{q+2H-2}\\
&\sim&s^{2H}2 H(2H-1)\sum_{q=0}^{2m+3} \frac{d_q}{q+2H-1}.
\end{eqnarray*}

\section{Application \label{sec4} }

\subsection{Bias for scaling signals }
In \citet{Kantelhardt:2001dr} the bias (of the DFA fluctuation function) for Hurst  exponents $H=0.5,0.65,0.9$ was found by means of Monte Carlo simulations. From this bias they proposed the modified DFA fluctuation function
\begin{equation}
F^2_\textrm{mod}(s)=\frac{F^2(s)}{K^2 (s)},
\label{modF}
\end{equation}
with 
\begin{equation*} 
K^2(s)= \frac{\mathbb E F^2 (s) \tau^{2H} }{\mathbb E F^2 (\tau) s^{2H}}.
\end{equation*}
If we assume $\tau$ is large, such that $\mathbb E F^2 (\tau)=\lambda_{m,H} \tau^{2H}$ holds (approximately), then 
\begin{equation} 
K^2(s)= \frac{\mathbb E F^2 (s)  }{\lambda_{m,H} s^{2H}},
\label{KF}
\end{equation}
which implies $\E F^2_\textrm{mod}(s)= \lambda_{m,H}s^{2H}$. 

\begin{figure}
\begin{center}
\includegraphics[width=\linewidth]{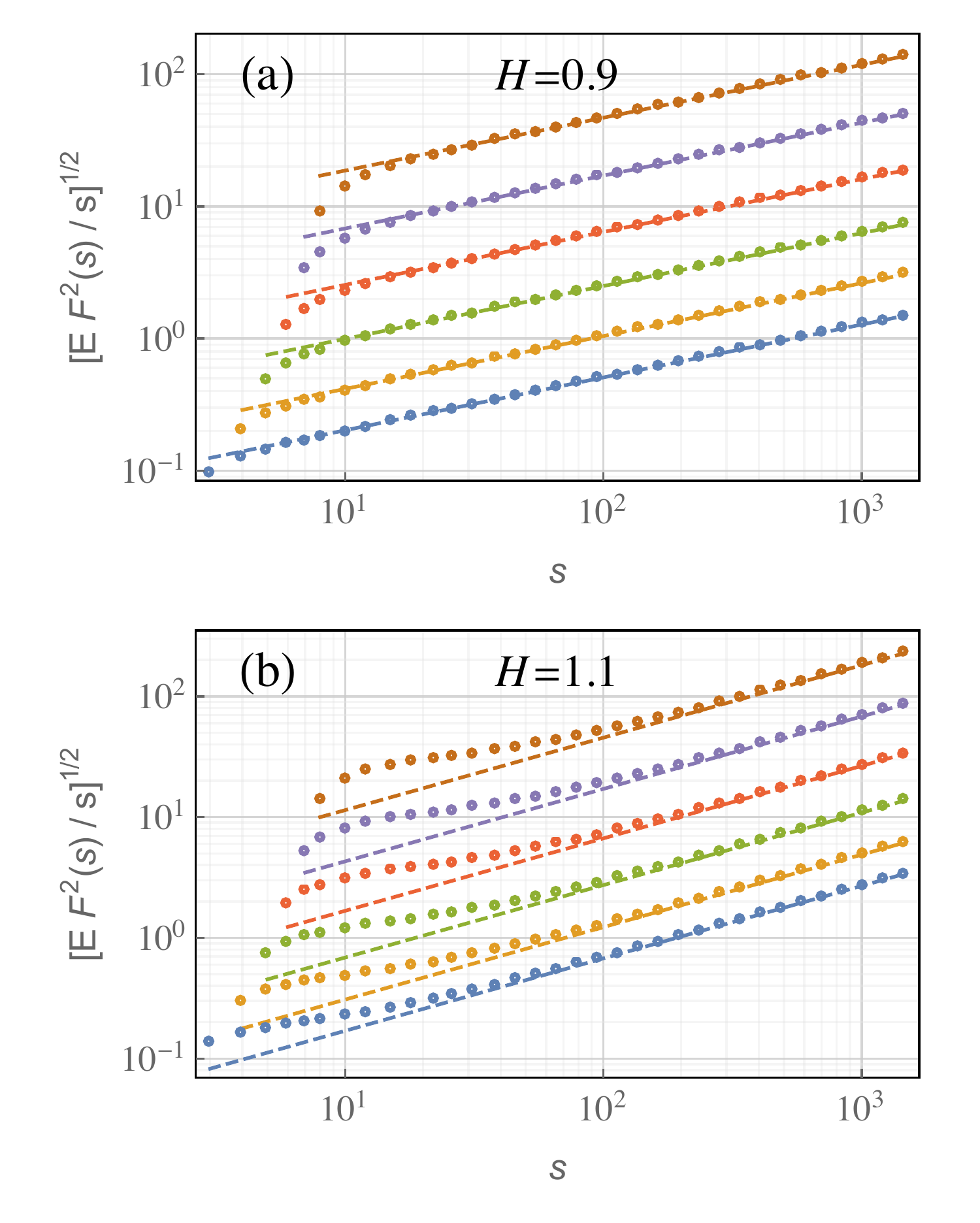}
\caption{Detrended fluctuation analysis for input signals with Hurst exponent (a) $H=0.9$ and (b) $H=1.1$. 
In both figures (a,b) the graphs from bottom to top corresponds to DFA of increasing order $m$, from $m=1$ (bottom) to $m=6$ (top). Dashed lines are the asymptotic scaling $\lambda_{m,H}^{1/2}s^H$ (see text).
The squared fluctuation functions have been shifted by factors $10^{m-1}$.
}
\label{figA}
\end{center}
\end{figure}

Eq.~\eqref{EACF} can be used  to compute the bias, i.e., the difference between the fluctuation function and asymptotic scaling, for signals with Hurst exponents  $0<H<1$. An example is shown in Fig.~\ref{figA}a where we have used $H=0.9$. This gives similar result as provided by \citet{Kantelhardt:2001dr} (see their Fig. 2a) since the only difference between the analytical and Monte Carlo method is the negligible error caused by finite sample length for the latter.

  Eq.~\eqref{ES}  can also be used to compute the bias for signals with Hurst exponents $1<H<2$.  Fluctuation functions with corresonding asymptotic scaling for $H=1.1$ is shown in Fig.~\ref{figA}b. 

\begin{figure}
\begin{center}
\includegraphics[width=\linewidth]{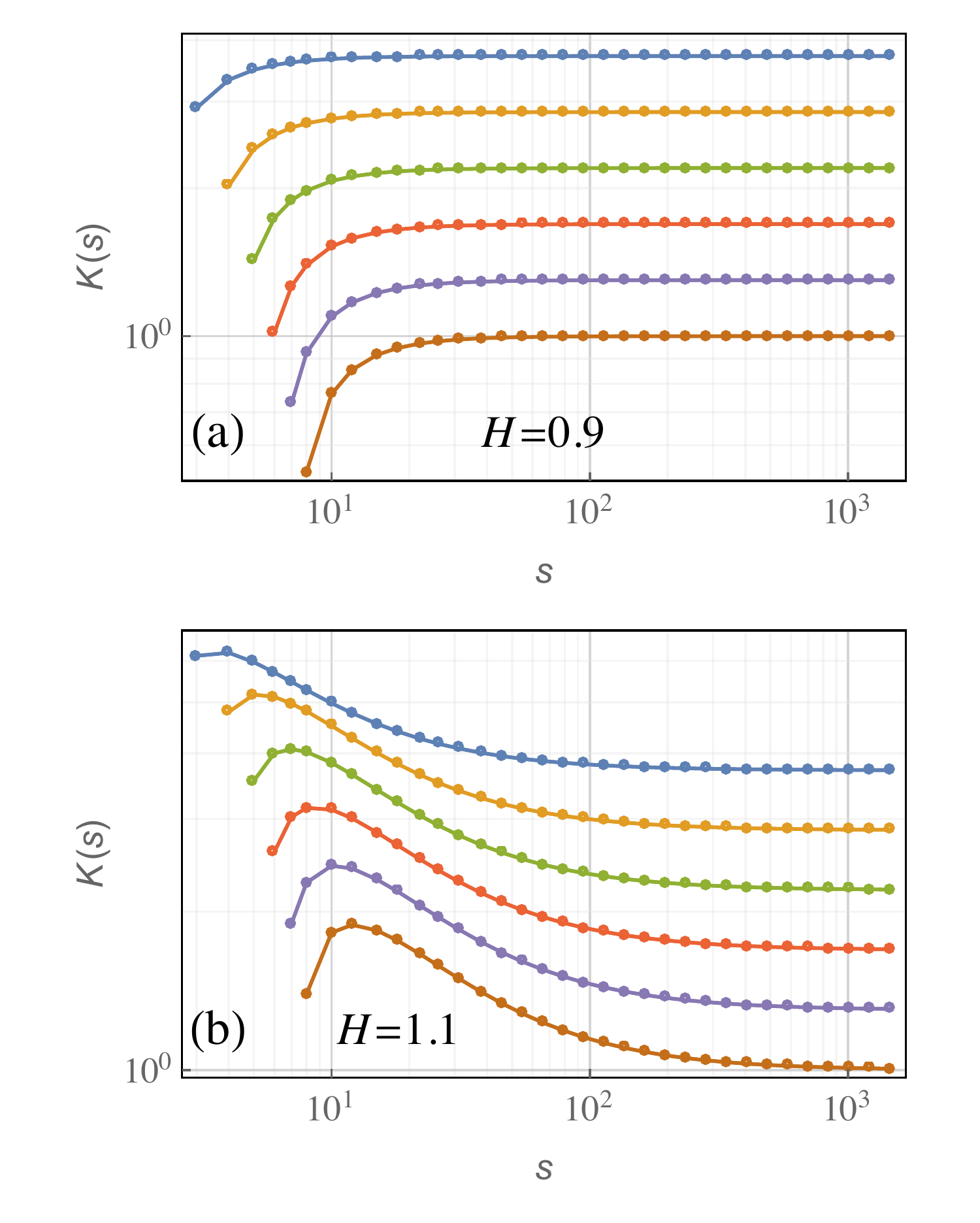}
\caption{ The correction function $K(s)$ for input signals with Hurst exponent (a) $H=0.9$ and (b) $H=1.1$
In both figures (a,b) the graphs from bottom to top corresponds to DFA of decreasing order $m$, from $m=6$ (bottom) to $m=1$ (top).
The correction functions has been shifted by factors $1.3^{6-m}$.
}
\label{figB}
\end{center}
\end{figure}

The correction functions Eq.~\eqref{KF} for $H=0.9$ and $H=1.1$ are shown in Fig.~\ref{figB}. A practical problem is that $K(s)$ depends on the (unknown) Hurst exponent. In \citep{Kantelhardt:2001dr} this dependence was found to be weak for $H=0.5,0.65,0.9$. Based on this finding the authors suggested to use Eq.~\eqref{modF} with the correction function for $H=0.5$. While using this modified DFA will improve the scaling for $H<1$, it will actually increase the bias for signals with Hurst exponents $H>1$. For $H=0.9$ and $H=1.1$, this can be seen from Figs.~\ref{figA} and \ref{figB}, where we observe that the bias has different signs. 



\subsection{Ornstein Uhlenbeck processes}
\begin{figure}
\begin{center}
\includegraphics[width=\linewidth]{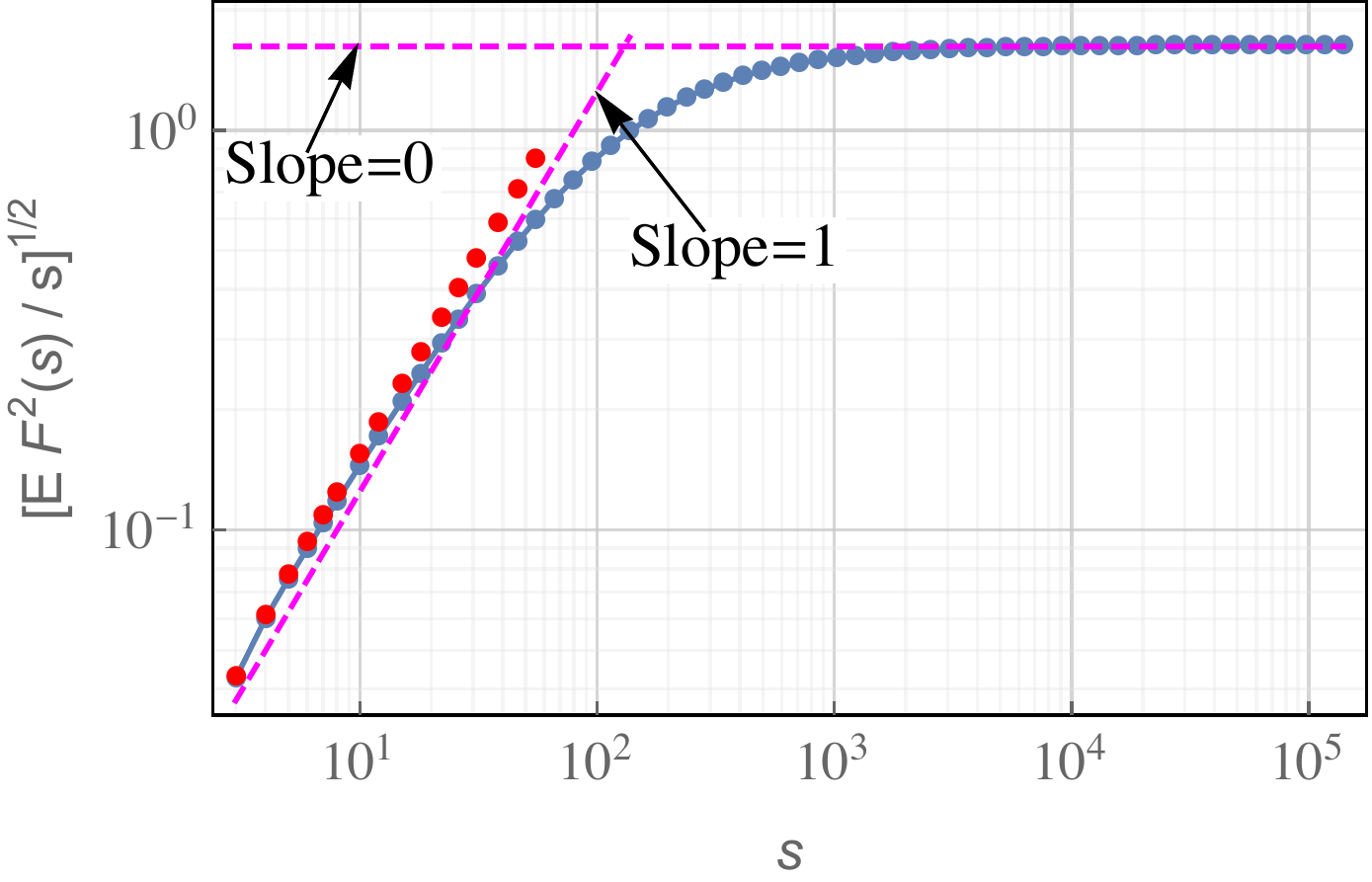}
\caption{Detrended fluctuation analysis of order 1 for Ornstein Uhlenbeck process with characteristic correlation time $\tau=20$ (blue curve). The red points are DFA1 for random walks (Hurst exponent $H=1.5$). }
\label{figOU}
\end{center}
\end{figure}

Another application is to study the behaviour of DFA for signals that are not scaling. Here we consider the class of Ornstein-Uhlenbeck (OU) processes. An OU is the solution to the Langevin equation
\begin{equation}
dX(t) = -\frac{1}{\tau} X(t)dt +\sigma dB(t),
\label{OU}
\end{equation}
where $B(t)$ is a standard ($\mathbb E B(1)^2=1$) Brownian motion, $\sigma>0$ is a scale parameter and $\tau>0$ is the characteristic correlation time. 
We choose initial condition such that $X(t)$ is stationary. This implies that the auto-covariance function takes the form
\begin{equation}
\E X(t)X(s)=\frac{\exp(-|t-s|/\tau)}{2 \sigma^2}
\label{CovOU}
\end{equation}
Again, we can use Eq.~\eqref{EACF} to  compute the expected value of the squared DFA fluctuation. An example is shown in Fig.~\ref{figOU}.

While OU processes do not have well-defined Hurst exponents as defined in Section~\ref{sec2}, the second moment scales asymptotically: On long time scales ($\tau \rightarrow \infty$), $X(t)$ is white noise, while on short time scales ($\tau \rightarrow 0$), $X(t)$ converges to a Brownian motion. Thus, for the DFA fluctuation function we should expect a scaling exponent close to $H=0.5$ on long time scales. It is seen in Fig.~\ref{figOU} that this holds. On small time scales there exists bias in DFA for signals that exhibit scaling behavior.  Relevant here is the bias for random walks ($H=1.5$).  In Fig.~\ref{figOU} it is seen that the OU DFA fluctuation function, with $\tau=20$, is consistent with random walks on small time scales. 

We note that an AR(1) is an discretised OU process, and more results on the AR(1) DFA fluctuation function can be found in \citep{Holl:2015gv}.

\subsection{Missing data \label{secmd}}
Based on Eq.~\eqref{F3} we can modify DFA to handle missing data. Define $\delta(t)$ to be  zero if $X(t)$ is missing and unity otherwise. We make the assumption that at least one $X(t+k)X(t+j)$ is non-missing. A sufficient, but not necessary, condition for this to hold is that at least one window contain no gaps.   
Let 
 $$p_{k,j}= \frac{\#\,   \text{of windows }}{\#\,   \text{ of non-missing }  X(t+k)X(t+j)}. $$
We propose the estimator
\begin{eqnarray}
\hat F^2_t(s)&=&- \frac{1}{2s}\sum_{k=1,j=1}^s p_{k,j} a_{k,j} \nonumber\times\\
&& [X(t+k)-X(t+j)]^2\times \nonumber \\
&& \delta(t+k) \delta(t+j).
\label{Fhat}
\end{eqnarray}
We define $\hat F^2(s)$ to be the average of $\hat F^2_t(s)$ (averaging over the different windows $t$ used). 
Without missing data the \textit{modified} fluctuation function $\hat F(s)$ is the same as the fluctuation function $F(s)$ in the gap-free case. 
For a time series with gaps $\hat F(s)$ is equal in expectation to $F(s)$.

The equality $\E \hat F^2(s) =\E  F^2(s)$ holds if the input signal is stationary or non-stationary with stationary increments:
Applying the expectation operator on Eq.~\eqref{Fhat} we have
\begin{eqnarray*}
\E \hat F^2_t(s)&=&- \frac{1}{2s}\sum_{k=1,j=1}^s p_{k,j} a_{k,j} S(|k-j|)\times\\
&& \delta(t+k) \delta(t+j),
\end{eqnarray*}
and since at least one $\delta(t+k) \delta(t+j)$ is assumed non-zero, the equality $\E \hat F^2(s) =\E  F^2(s)$ follows.  Whereas the ordinary fluctuation function is always non-negative, the modified fluctuation function can become negative. Practically, one can resolve this problem by letting the fluctuation function be undefined if negative values occur.

\begin{figure}
\begin{center}
\includegraphics[width=\linewidth]{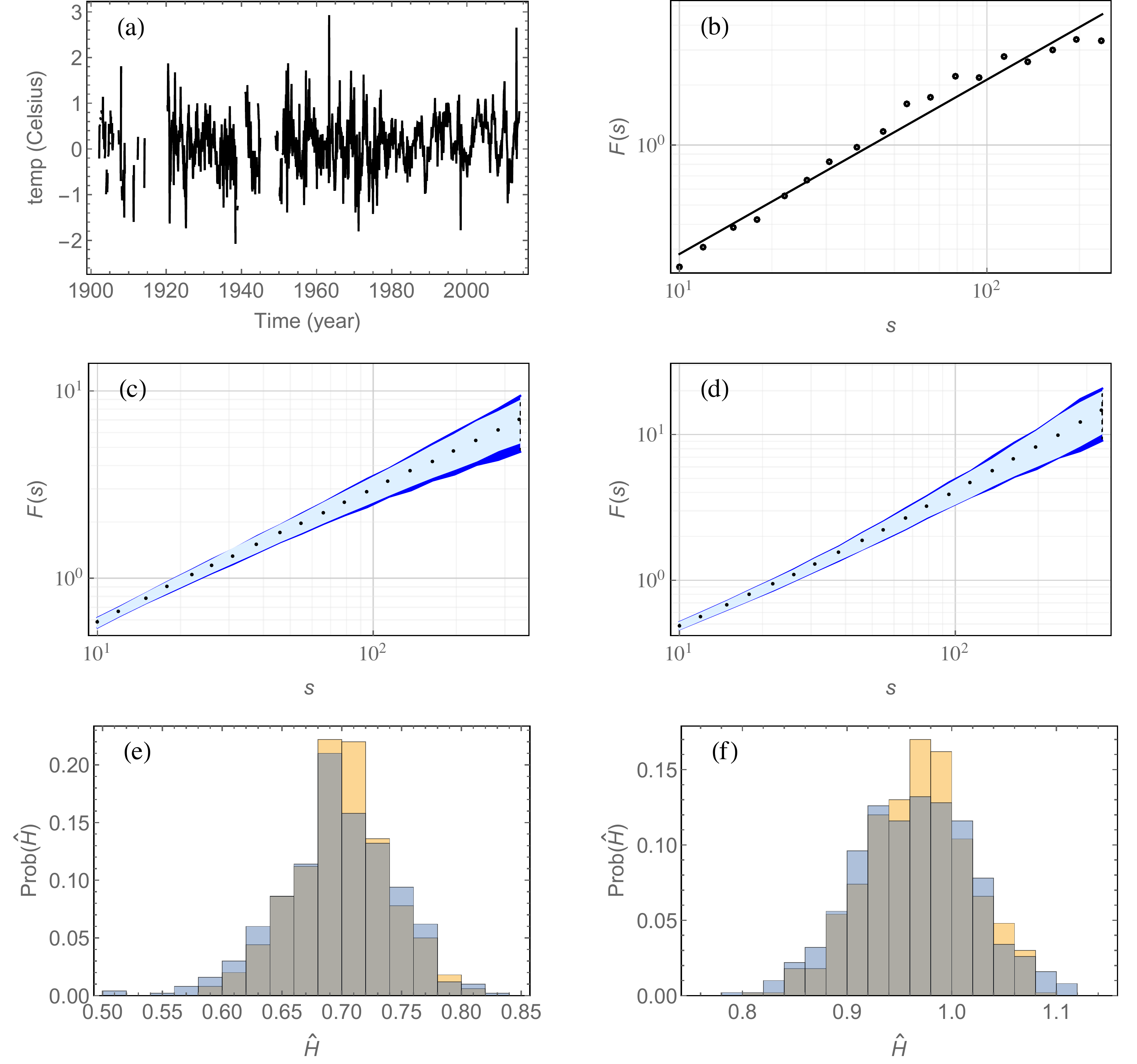}
\caption{a) Monthly reconstructed temperature for the $5^\circ \times 5^\circ$ grid centered at 7.5$^\circ$S, 172.5$^\circ$W. b) Modified Detrended Fluctuation Analysis of the time series in (a). (c-f) The results of the Monte Carlo study of DFA explained in section \ref{secmd}. In (c-d) the black points are the ensemble means of the modified and ordinary fluctuation functions, light blue is the $90\%$ confidence interval for the ordinary fluctuation function, and light plus dark blue is the $90\%$ confidence interval for the modified fluctuation function. (e-f) Blue plus gray is the probability density function (PDF) of ensembles of estimated Hurst exponents according to the modified DFA, while yellow plus gray is the PDF according to the ordinary DFA. 
}
\label{figMissing}
\end{center}
\end{figure}

Examples of time series with gaps are some of the regional temperatures  analyzed in \citet{Lovsletten:2016bja}. We use one of these time series, from the HADCRUT4 data product \citep{Morice:2012dw}, to demonstrate the usage of the modified DFA to handle missing data. The chosen time series is the surface temperature in the tropical  Pacific, and is shown in Fig.~\ref{figMissing}a. The modified DFA of this series is presented in Fig.~\ref{figMissing}b. In this work we use non-overlapping windows starting from the left, and (modified) DFA of order $m=2$.
 We observe that the fluctuation function is rather poorly approximated by a power-law. This is an expected result since in \citep{Lovsletten:2016bja} we showed that an AR(1) model is significantly better than an fGn (power-law) model for this temperature series due to the influence of the El Ni{\~ n}o Southern Oscillation (see discussion in \citep{Lovsletten:2016bja}). 

Using the same gap-sequence as in Fig.~\ref{figMissing}a, we compare the modified DFA with the ordinary DFA by computing the fluctuation function from an ensemble of 500 fGn's with Hurst exponent $H=0.7$ and sample size $n=1368$ (same sample size as in Fig.~\ref{figMissing}a).  For each ensemble member, data points are omitted to construct the same gap sequence as in Fig.~\ref{figMissing}a and then  the modified DFA fluctuation function is computed. In Fig.~\ref{figMissing}c the results are presented in form of ensemble means and $90\%$ confidence intervals. The ensemble means are visually indistinguishable for the modified and ordinary fluctuation functions, while the uncertainty of the former is increased due to the gaps. 

The increased uncertainty is also seen in the estimated Hurst exponents, see  Fig.~\ref{figMissing}e. Here, the estimators are the slopes from linear regression of $\log F(s)$ (and $\log \hat F(s)$) against $\log s$.

The same Monte Carlo experiment,  using an ensemble of fBm's with Hurst exponent $H=1.1$, yields similar results which are presented in Fig.~\ref{figMissing}d and f. 

An alternative estimator, based on Eq.~\eqref{F2}, is 
\begin{eqnarray}
\tilde F^2_t(s)&=&\frac{1}{s}\sum_{k=1,j=1}^s p_{k,j} a_{k,j}\times \nonumber\\
&& X(t+k)X(t+j)\times \nonumber \\
&& \delta(t+k) \delta(t+j),
\label{Ftilde}
\end{eqnarray}
and $\tilde F^2(s)$ defined as the average of $\tilde F^2_t(s)$. It is straight-forward to verify that $\E \tilde F^2(s) =\E  F^2(s)$ for stationary input signals. 
However, for input signals that are non-stationary with stationary increments,  $\tilde F(s)$ does not have the desirable property of equality in expectation to the fluctuation function $F(s)$ in the gap-free case. As an example, consider a input signal with Hurst exponent $1<H<2$. In the gap-free case, the time-dependent part of the expected squared fluctuation vanishes, i.e., 
\begin{equation*}
\mathbb E F_t^2(s)= \frac{1}{s}\sum_{k=1}^s\sum_{j=1}^s a_{k,j} (|t+k|^{2h}+|t+j|^{2h})=0,
\end{equation*}
(see Eq.~\eqref{null} in the appendix). For $\tilde F^2_t(s)$ the time-dependent part will not, in general, vanish. This is due to the additional multiplicative factors $p_{k,j}\delta(t+k) \delta(t+j)$ in Eq.~\eqref{Ftilde}.

%


\section{Concluding remarks}
In this paper, several new propositions for DFA have been formulated and proven. These include the relationship between the DFA fluctuation function and the acvf and variogram, derived from the sample forms Eqs.~\eqref{F2} and \eqref{F3}.  The results were derived under the assumption that the input signal in DFA is either  stationary or non-stationary with stationary increments, or one of these  superposed on a polynomial trend of order less than the order of DFA (results on trends not accounted for by DFA can be found in \citet{Hu:2001ia} and \citet{Kantelhardt:2001dr}). For these classes of input signals  the present paper has established  that the residual variance in different windows are equal in expectation.
The power-law scaling of the DFA fluctuation function has been rigorously  proven for stochastic processes with Hurst exponents $H\in\{(0,1)\cup (1,2) \}$.

It has also been demonstrated for these classes of signals that the new developments of the DFA method can be used to compute analytically the bias of the DFA-estimate. For this purpose we used the weight functions and asymptotic weight functions. The Mathematica code for these functions are found as supplementary material to this article. Therein, it is demonstrated that the modified fluctuation function proposed in \citet{Kantelhardt:2001dr} degrades the scaling property for scaling input signals with Hurst exponent greater  than unity.

From an applied physics point of view, the most useful result of this study may be the  method of handling missing data in DFA proposed in Eq.~\eqref{Fhat}. For ensemble averages, the modified DFA fluctuation function  with missing data is the same as the DFA fluctuation function without missing data.

Some of the theory presented in this paper is probably a suitable starting point to prove the correctness of the multfractal DFA introduced by \citet{Kantelhardt:2002cc}, as well as the variance and limiting distribution of the DFA fluctuation function. 

%
%

\acknowledgments This work has received support from the Norwegian Research Council under Contract 229754/E10. The author thanks K. Rypdal and M. Rypdal for useful discussions and language editing.

\appendix

\section{Simple proofs \label{VA}}
Recall the definition of the weight matrix  
$$A=D^T(I-Q)D$$
and hat matrix
$$Q=B^T\left(BB^T\right)^{-1} B. $$
Since $Q$ is a projection matrix, vectors $\mathbf v$ that are in the row-space of $B$ will be mapped to itself, i.e. $Q\mathbf v=\mathbf v$, and thus $(I-Q)\mathbf v =0$. 

\eqref{DFAdef} $\Leftrightarrow$ \eqref{F2}:
Let $\mathbf 1$ be a $(s\times 1)$ vector of ones.
For $t>0:$ 
$$\mathbf Y(t)= D\mathbf{X}(t)+ \mathbf 1  Y(t).$$
The proof is completed by noting that $\mathbf 1$ is in the the row-space of $B$.\\~\\
\eqref{F2}$\Leftrightarrow$\eqref{F3}:
This equality holds if 
\begin{equation}
\sum_{k=1,j=1}^s a_{k,j} (X(t+k)^2+X(t+j)^2)=0 
\label{null}
\end{equation}
Fix $k$ and consider the sum 
$$\sum_{j=1}^s a_{k,j}X(t+k)^2,$$
which is element $k$ of the vector  $A\mathbf 1 X(t+k)^2$. We have $D\mathbf 1= (1,2,\ldots,n)^T$, which is in the row-space of $B$. Thus $A\mathbf 1 X(t+k)^2=0$. Since this holds for all $k=1,2,\ldots s$, we can conclude that 
$$\sum_{k=1,j=1}^s a_{k,j}X(t+k)^2=0.$$
Since $A$ is a symmetric matrix we also have
$$\sum_{k=1,j=1}^s a_{k,j}X(t+j)^2=0.$$

In the sequel it is shown that the relationships between the DFA fluctuation function and the acvf and variogram estimators, Eqs.~\eqref{acvf1} and \eqref{ES}, and the power-law scaling of the DFA fluctuation function Eq.~\eqref{MR2}, remain valid when a polynomial  trend is superposed on the signal. 
Let $Z(t)$ be a stochastic process with mean zero and acvf $\gamma(t,s)$. It is assumed that $Z(t)$ is either stationary or non-stationary with stationary increments. Define
$$T(t)=\beta_0+\beta_1 t+\ldots+\beta_q t^q,\quad t=1,\ldots,n,$$
where $q$ is an integer in the range $0 \leq q\leq  m-1$.  
Let
$$X(t)=T(t)+Z(t).$$
By Eq.\eqref{F2} we have
\begin{eqnarray}
F^2_t(s)&=&\frac{1}{s}\sum_{k=1,j=1}^s a_{k,j} Z(t+k)Z(t+j)\nonumber\\
&+&\frac{1}{s}\sum_{k=1,j=1}^s a_{k,j} Z(t+k)T(t+j)\nonumber \\
&+&\frac{1}{s}\sum_{k=1,j=1}^s a_{k,j} T(t+k)Z(t+j)\nonumber\\
&+&\frac{1}{s}\sum_{k=1,j=1}^s a_{k,j} T(t+k)T(t+j) \label{Ftrend}.
\end{eqnarray}
Since $\E Z(t)=0$, the middle terms vanish in expectation, and thus
\begin{eqnarray}
\E F^2_t(s)&=&\frac{1}{s}\sum_{k=1,j=1}^s a_{k,j} \gamma(t+k,t+j)\nonumber\\
&+&\frac{1}{s}\sum_{k=1,j=1}^s a_{k,j} T(t+k)T(t+j). 
\end{eqnarray}
We have
\begin{eqnarray}
\sum_{k=1,j=1}^s a_{k,j} T(t+k)T(t+j) =\mathbf{T}(t)^T A \mathbf{T}(t), 
\end{eqnarray}
where $\mathbf{T}(t)=[T(t+1),\ldots,T(t+s)]^T$. One can use the formulas for sums of powers, e.g. \cite{W}, to verify that $D\mathbf{T}(t)$ is in the row-space of $B$. Hence 
\begin{eqnarray*}
\sum_{k=1,j=1}^s a_{k,j} T(t+k)T(t+j) =0,
\end{eqnarray*}
and thus 
\begin{eqnarray}
\E F^2_t(s)&=&\frac{1}{s}\sum_{k=1,j=1}^s a_{k,j} \gamma(t+k,t+j).\nonumber
\end{eqnarray}

\section{Asymptotic weight function \label{VB}}
The weight matrix can be written
$$A=D^TD-D^TQD,$$
where $Q$ is the hat matrix defined in Eq.~\eqref{hatmatrix}. Element $(i,j)$ of $D$ is one if $i\geq j$ and zero otherwise. Thus, element $(i,j)$ in the first matrix product is 
$$(D^TD)_{i,j}=s+1-\max\{i,j\}.$$
Summing the $j$'th (sub)-diagonal yield
\begin{eqnarray}
\sum_{k=1}^{s-j}(D^TD)_{k,k+j}& = & \sum_{k=1}^{s-j}(s+1-(k+j))\nonumber \\
&=&s^2/2-s j+s/2   +j^2/2-j/2\nonumber\\
&\sim&s^2/2-s j+j^2/2. \label{DD}
\end{eqnarray}
Computing  the term $D^TQD$ is more tedious, but straight-forward.  The starting point is  the hat matrix $Q$. 
Denote by $(B B^T)^{-1}_{i,j}$ element $(i,j)$ of the inverse of $B B^T$. By observing that column $j$ of $B$ is $(j^0,j^1,\ldots,j^m)$, it is seen  that
$$Q_{p,q}= \sum_{d=1,l=1}^{m+1}  p^{d-1}q^{l-1}(B B^T)^{-1}_{d,l}.$$

\begin{eqnarray}
[D^TQ D]_{i_1,i_2}= \sum_ {k_1=i_1}^s \sum_ {k_2=i_2}^s Q_{k_1,k_2}\nonumber&&\\
=\sum_ {k_1=i_1}^s \sum_ {k_2=i_2}^s\sum_{d=1,l=1}^{m+1}  k_1^{d-1}k_2^{l-1}(B B^T)^{-1}_{d,l}\nonumber\\
\sim \sum_{d=1,l=1}^{m+1}  (s^{d}-i_1^{d})  (s^{l}-i_2^l)/(dl)(B B^T)^{-1}_{d,l}
\label{DHD}
\end{eqnarray}

Using the asymptotic expression of $B B^T$,
$$ (B B^T)_{i,j} =\sum_{t=1}^s t^{i+j-2}\sim \frac{s^{i+j-1}}{i+j-1},$$
one can use the definition of the inverse matrix to verify that
\begin{equation}
(B B^T)^{-1}_{d,l}\sim \tilde c_{d,l}/s^{d+l-1}
\label{asymBinv}
\end{equation}
Inserting Eq.~\eqref{asymBinv} into Eq.~\eqref{DHD} yields

$$ [D^TQ D]_{i_1,i_2}\sim 
\sum_{d=1,l=1}^{m+1}  (s^{d}-i_1^{d})  (s^{l}-i_2^l)c_{d,l}/s^{d+l-1},$$
where we have defined   $c_{d,l}=\tilde c_{d,l}/(dl)$. Summing the $j$'th (sub)-diagonal yield

\begin{eqnarray}
\sum_{k=1}^{s-j}(D^TQD)_{k,k+j}&\sim& \nonumber\\
\sum_{k=1}^{s-j} \sum_{d=1,l=1}^{m+1}  (s^{d}-k^{d})  (s^{l}-(k+j)^l)c_{d,l}/s^{d+l-1}& \sim & \nonumber\\
\sum_{d=1,l=1}^{m+1} c_{d,l}\left (s^2-sj -\frac{s^2}{l+1}+\frac{s^{1-l} j^{l+1}}{l+1}\right)\label{b1}\\
-\sum_{d=1,l=1}^{m+1} c_{d,l}\frac{s^{-d+1}(s-j)^{d+1} }{d+1}  \\
 +\sum_{d=1,l=1}^{m+1} c_{d,l}\sum_{r=0}^l \binom{l}{r}\frac{(s-j)^{d+l+1-r} j^r}{(d+l+1-r)s^{d+l-1}}\label{b3}\\
 = \sum_{q=0}^{2m+3}s^{2-q}j^q b_q.\label{b5}
\end{eqnarray}
The terms \eqref{b1}-\eqref{b3} can be written $$\sum_{q=0}^{2m+3 }s^{2-q}j^q b^{(k)}_q,\quad, k=1,2,3,$$ respectively. This implies the equality Eq. \eqref{b5}, with 
$$ b_q=b^{(1)}_q+b^{(2)}_q+b^{(3)}_q.$$
The coefficients  $b^{(k)}_q$, found by re-organising terms, are given by:
$$b^{(1)}_q
=
\begin{cases}
\sum_{d=1,l=1}^{m+1}c_{d,l}-c_{d,l} \frac{1}{l+1}& \text{if     } q=0, \\
 -\sum_{d=1,l=1}^{m+1}c_{d,l}& \text{if     } q=1, \\
\frac{1}{q}\sum_{d=1}^{m+1}c_{d,q-1}& \text{if     } 2\leq q \leq m+2,\\
0& \text{if     } q>m+2.
\end{cases}
$$

$$b^{(2)}_q
=
\begin{cases}
-\sum_{d=1,l=1}^{m+1}\frac{c_{d,l}}{d+1}& \text{if     } q=0, \\
\sum_{d=1,l=1}^{m+1}\frac{c_{d,l}}{d+1}\binom{d+1}{d} & \text{if     } q=1, \\
(-1)^{q-1}\sum_{d=q-1,l=1}^{m+1}\frac{c_{d,l}}{d+1}\binom{d+1}{d+1-q} & \text{if     } 2\leq q \leq m+2,\\
0& \text{if     } q>m+2.
\end{cases}
$$

$$b^{(3)}_q=\sum_{\substack{d+l\geq q-1,\\d\geq1, l\geq 1}}^{m+1} a_q^{(d,l)} c_{d,l}, $$

$$ a_k^{(d,l)}= \sum_{r=0}^{\min{\{l,k\}}} \binom{l}{r} \frac{(-1)^{k-r}}{d+l+1-r} \binom{d+l+1-r}{d+l+1-k}. $$
Using Eq.~\eqref{DD} and \eqref{b5}, the coefficients follows:
$$d_q
=
\begin{cases}
1/2-b_0& \text{if     } q=0, \\
-1-b_1& \text{if     } q=1, \\
1/2-b_2 & \text{if     }  q =2,\\
-b_q& \text{if     } q>2.
\end{cases}
$$
\vspace{1cm}


%

\end{document}